\documentclass[preprint,12pt]{elsarticle}
\usepackage{graphicx} % Required for inserting images

\usepackage[utf8]{inputenc}
\usepackage[english]{babel}
\usepackage{graphicx}
\usepackage{amssymb}
\usepackage{amsmath}
\usepackage{tikz}
\usetikzlibrary{decorations.pathreplacing} %for the braces
\usepackage{commath}
\usepackage{upgreek}
\usepackage[hidelinks]{hyperref}
\usepackage{color}
\usepackage{wrapfig}
\usepackage{graphicx,wrapfig,lipsum}
\graphicspath{ {./images/} }

\usepackage{amsthm}
\newtheorem{theorem}{Theorem}[section]

\newtheorem{lemma}[theorem]{Lemma}

\date{\vspace{-5ex}}

\usepackage[a4paper, total={6 in, 9.5 in}]{geometry}

\journal{Nonlinear Analysis}

\begin{document}

\begin{frontmatter}

\title{Asymptotics and periodic dynamics in a negative chemotaxis system with cell lethality %
% titulo puesto en el otro articulo como preprint
%Periodicity in a fully parabolic negative chemotaxis system with source term and lethal interaction
}

\author[ucm,upb,icai]{Federico Herrero-Hervás \corref{cor1}}
\ead{fedher01@ucm.es}
\author[ucm]{Mihaela Negreanu}
\cortext[cor1]{Corresponding author}

\affiliation[ucm]{organization={Instituto de Matemática Interdisciplinar, Departamento de Análisis Matemático y Matemática Aplicada, Universidad Complutense de Madrid},
            city={Madrid},
            postcode={28040},
            country={Spain}}

\affiliation[upb]{organization={Universität Paderborn},
            city={Paderborn},
            postcode={33098},
            country={Germany}}

\affiliation[icai]{organization={Departamento de Matemática Aplicada, ICAI, Universidad Pontificia de Comillas},
            city={Madrid},
            postcode={28015},
            country={Spain}}            
\begin{abstract}
This work studies the following system of parabolic partial differential equations
\begin{equation*}
\begin{cases}
       \displaystyle \frac{\partial u}{\partial t} = D\Delta u + \chi \nabla \cdot(u \nabla v) + ru(1-u) - u v, \quad & x \in \Omega, ~t > 0, \\\\
      \displaystyle \frac{\partial v}{\partial t} = \Delta v + a u -v+ f(x,t),  \quad & x \in \Omega, ~t > 0,
\end{cases}
\end{equation*}
modeling the negative chemotaxis interactions between a biological species and a lethal chemical substance that is supplied according to the known function $f(x,t)$. 
\\\\
It is shown that if $f$ converges to a spatially homogeneous function $\tilde{f}$ in a certain sense, then the solution $(u,v)$ satisfies
$$
||u-\tilde{u}||_{L^2(\Omega)} + ||v-\tilde{v}||_{L^2(\Omega)} \to 0 \quad \text{as } t \to \infty,
$$
where $(\tilde{u},\tilde{v})$ is the solution to the associated ODE system
\begin{equation*}
    \begin{cases}
  \displaystyle  \frac{d \tilde{u}}{dt~} = r \tilde{u} (1 - \tilde{u}) - \tilde{u}\tilde{v}, \quad & t>0,\\\\
  \displaystyle  \frac{d \tilde{v}}{dt~} = a\tilde{u} - \tilde{v} + \tilde{f},\quad & t>0.
\end{cases}
\end{equation*}
Some final remarks are given for the case in which $\tilde{f}$ is a time periodic function, and under which hypotheses do $(\tilde{u},\tilde{v})$ inherit this periodicity.

\end{abstract}
\end{frontmatter}

\section{Introduction}
This work is centered around the asymptotic behavior of the solution to the following system of two coupled nonlinear parabolic partial differential equations
\begin{equation} \label{1.1}
\begin{cases}
       \displaystyle \frac{\partial u}{\partial t} = D\Delta u + \chi \nabla \cdot(u \nabla v) + ru(1-u) - u v, \quad & x \in \Omega, ~t > 0, \\\\
      \displaystyle \frac{\partial v}{\partial t} = \Delta v + a u -v+ f(x,t),  \quad & x \in \Omega, ~t > 0,
\end{cases}
\end{equation}
subject to nonnegative initial values and Neumann homogeneous boundary conditions
\begin{equation}\label{1.1.2}
    \begin{cases}
        u(x,0) = u_0(x) \geq 0, \quad v(x,0) = v_0(x) \geq 0, & x \in \Omega,\\\\
        \displaystyle \frac{\partial u}{\partial \nu} = \frac{\partial v}{\partial \nu} = 0, & x \in \partial \Omega, ~t>0,
    \end{cases}
\end{equation}
posed over a bounded and convex domain $\Omega \subset \mathbb{R}^n$ with smooth boundary, with parameters $D, ~\chi, ~ r >0,~a\geq0$ and a nonnegative known source function $f \geq 0$.
\\\\
System \eqref{1.1} has recently been proposed in~\cite{HN25} as a model for the interaction between a bacterial population and a lethal chemical agent, such as a toxin or an antibiotic. In the model, $u(x,t)$ denotes the bacterial density and $v(x,t)$ the concentration of the chemical at position $x \in \Omega$ and time $t>0$, accounting for the following mechanisms:
\begin{enumerate}
    \item \textit{Cell motility and chemical diffusion}, incorporated by the Laplacian terms $\Delta u$ and $\Delta v$.
    \item \textit{Negative chemotaxis}, this is, bacterial repulsion from regions of high chemical concentration, described by $\chi \nabla \cdot (u\nabla v)$.
    \item \textit{Logistic bacterial growth and chemical lethality}, modeled respectively by the terms $r u(1-u)$ and $-u v$.
    \item \textit{Self-production and degradation of the chemical}, captured by $a u - v$.
    \item \textit{An external supply of the substance}, represented by the source function $f$.
\end{enumerate}
These interactions considered for the model are common in microbiological systems where organisms react to toxins, antibiotics, or self-produced metabolic byproducts. This is the case for instance for \textit{E. coli} bacteria and $H_2O_2$. It is known that \textit{E. coli} cells naturally produce $H_2O_2$ molecules \cite{h1,h2}, which act as a lethal chemorepellent \cite{bnv,u,im}, leading to bacterial migration towards regions of lower concentrations in order to survive.
\\\\
Although for most bacteria such as \textit{E. coli}, the self production of a chemo-repellent is in general not a means of self regulation or signaling ---but instead a consequence of their metabolism \cite{nature}---, its combination with the negative chemotaxis and the intrinsic dynamics of self-production, degradation, logistic growth and lethality creates an intricate feedback mechanism between bacterial motion and the interactions with the chemical, whose asymptotics are analyzed in this work.
\\\\
We note as well that by allowing $a=0$, the model also accounts for the interactions with an exogenous antibiotic, this is, when the chemical is not produced by the bacteria themselves.
\\\\
From a mathematical point of view, the first models for chemotaxis were introduced by Keller and Segel \cite{KS1, KS2}, considering the nonlinear term $\nabla \cdot(u\nabla v)$. Over the years, the Keller-Segel model has been extensively studied, showing the rich behavior of its solutions, ranging from aggregation phenomena to finite-time blow-up (see for instance the surveys by Horstmann \cite{Horstmann} or Hillen and Painter \cite{HillenPainter2009}).
\\\\
In terms of global existence, it is well known that adding a logistic growth term to the model, resulting in a system of the form
\begin{equation}\label{log-chem}
    \begin{cases}
     \displaystyle   \frac{\partial u}{\partial t} = D_1 \Delta u - \chi \nabla \cdot (u \nabla v) + ru- \mu u^2, \quad & x \in \Omega, ~t>0, \\\\
     \displaystyle  \tau \frac{\partial v}{\partial t} = D_2 \Delta v + u -\alpha v, & x \in \Omega, ~t>0,
    \end{cases}
\end{equation}
under Neumann homogeneous boundary conditions and nonnegative initial values, can suppress the finite-time blow-up of solutions that certain solutions to the minimal Keller-Segel system experiment. 
\\\\
In particular, in Tello and Winkler \cite{tello-winkler-2007}, the parabolic-elliptic case (corresponding to $\tau = 0$) is analyzed for \textit{positive} chemotaxis (i.e. attractant signaling, with $\chi >0$). The authors prove that system \eqref{log-chem} admits a unique solution that remains globally bounded for any parameter set, if the problem is posed over a bounded domain in dimension $n \leq 2$. Moreover, in the higher dimensional case, global boundedness still holds provided that the damping effect of the logistic term is sufficiently strong, specifically when $\mu > \displaystyle \frac{n-2}{2} \chi$, so that the quadratic death term $-\mu u^2$ dominates the aggregation induced by chemotaxis.
\\\\
For the fully parabolic version (with $\tau >0$), system \eqref{log-chem} was studied by Winkler in \cite{Winkler10} for a general $\chi \in \mathbb{R}$, thus accounting for \textit{positive and negative} chemotaxis, obtaining global boundedness of solutions in arbitrary dimensions for convex domains if $\mu$ is large enough. For the \textit{positive} chemotaxis case there exists $M>0$ such that if $\displaystyle \frac{\mu}{\chi}>M$, the unique non-trivial spatially homogeneous steady state of the system becomes globally asymptotically stable, and any solution $(u,v)$ converges to it, as proven in \cite{Winkler14}. 
\\\\
Our system \eqref{1.1} mainly differs from \eqref{log-chem} not only in considering the negative chemotaxis scenario, but also in the asymmetrical death term $-uv$, ---which does not have a counterpart on the $v$ equation, as cell death does not result in an increase of the chemical concentration---, and in the external supply $f$, which greatly affects the dynamics of the system.
\\\\
Several other works have addressed similar competitive terms, mainly in the form of Lotka-Volterra interactions between two species that are attracted by the same chemical. The elliptic case for the equation of the chemical was studied in \cite{tello-winkler-12} in a weakly competitive scenario, obtaining as well the convergence of solutions to the non-trivial spatially homogeneous steady state, under appropriate parameter restrictions. The parabolic case is considered as well in \cite{bw1,mz,zh} with similar results for the asymptotic behavior of solutions.
\\\\
A system closely related to ours is given by
\begin{equation}\label{ant}
    \begin{cases}
        \displaystyle \frac{\partial u}{\partial t} = \Delta u - \chi \nabla( u \nabla v) + \mu u(1+f(x,t)-u), \quad & x \in \Omega, ~t>0, \\\\
     \displaystyle  \tau \frac{\partial v}{\partial t} = \Delta v + u - v, & x \in \Omega, ~t>0,
    \end{cases}
\end{equation}
where the logistic term satisfies $1+f(x,t) >0$, which can be interpreted as a spatio-temporal dependent carrying capacity for the biological species. In particular, the system is analyzed by Negreanu, Tello and Vargas in \cite{per1}, with the assumption that $f(x,t)$ converges to a spatially homogeneous and time periodic function $\tilde{f}$, modeling seasonal fluctuations in the maximum population load that the environment can support. A generalization of system \eqref{ant} with two species in competitive interaction and periodic logistic dynamics was introduced and studied in \cite{per2}
\\\\
In particular, for system \eqref{ant} in \cite{per1}, under the assumptions that
\begin{equation}\label{ant-restricciones}
    \begin{split}
        f \in L^\infty \big(\Omega \times(0, \infty) \big), 
        \\\\
        \int_0^\infty \int_\Omega |\nabla f|^2 dx \hspace{0.1 cm} dt < \infty,
        \\\\
         \int_0^\infty ||f-\tilde{f}||_{L^1(\Omega)} \hspace{0.1 cm} dt < \infty,
    \end{split}
\end{equation}
as well as 
\begin{equation}\label{ant-restricciones2}
    \mu > \frac{\chi^2}{16} \max \left\{ \int_\Omega u_0 \hspace{0.1 cm} dx, ~ \frac{1}{\mu} \left( 1 + ||f||_{L^\infty(\Omega \times (0,\infty))} \right) \right\},
\end{equation}
it is proven that any solution $(u,v)$ to \eqref{ant} verifies
$$
||u-\tilde{u}||_{L^2(\Omega)} + ||v-\tilde{v}||_{L^2(\Omega)} \to 0, \quad \text{as } t \to \infty,
$$
where $(\tilde{u},\tilde{v})$ is the solution to the ODE system associated to system \eqref{ant}, this is
\begin{equation}\label{ant-ode}
    \begin{cases}
        \displaystyle \frac{d \tilde{u}}{dt~} =  \mu \tilde{u}\big(1+\tilde{f}(t)-\tilde{u}\big), \quad & t>0, \\\\
     \displaystyle  \tau \frac{d \tilde{v}}{dt~}  = \tilde{u} - \tilde{v}, & t>0,
    \end{cases}
\end{equation}
with initial values 
$$
\tilde{u}(0) = \int_\Omega u_0(x)\hspace{0.1 cm} dx, \quad \tilde{v}(0)= \int_\Omega v_0(x) \hspace{0.1 cm} dx.
$$
\textbf{Main results and structure of the work}\\\\
We recall our system, given by
\begin{equation*}
\begin{cases}
       \displaystyle \frac{\partial u}{\partial t} = D\Delta u + \chi \nabla \cdot(u \nabla v) + ru(1-u) - u v, \quad & x \in \Omega, ~t > 0, \\\\
      \displaystyle \frac{\partial v}{\partial t} = \Delta v + a u -v+ f(x,t),  \quad & x \in \Omega, ~t > 0,
\end{cases}
\end{equation*}
under Neumann homogeneous boundary conditions and nonnegative initial values. Our aim is to develop a similar analysis as done for system \eqref{ant}, considering that $f(x,t)$ converges in time (in a certain sense specified below) to a spatially homogeneous function $\tilde{f}(t)$. 
\\\\In this way, the main result of the work is Theorem \ref{t-per}, that ensures that if certain parameter conditions are met, any solution $(u,v)$ converges in $L^2(\Omega)$ to $(\tilde{u},\tilde{v})$, the corresponding solution to the ODE system associated to our system \eqref{1.1}, given by
\begin{equation} \label{1.3}
\begin{cases}
  \displaystyle  \frac{d \tilde{u}}{dt~} = r \tilde{u} (1 - \tilde{u}) - \tilde{u}\tilde{v}, \quad & t>0,\\\\
  \displaystyle  \frac{d \tilde{v}}{dt~} = a\tilde{u} - \tilde{v} + \tilde{f},\quad & t>0,
\end{cases}
\end{equation}
with initial values
\begin{equation}\label{1.4}
    \tilde{u}(0) = \int_\Omega u_0(x) ~dx, \quad \tilde{v}(0) = \int_\Omega v_0(x) ~dx.
\end{equation}
This extends the existing information about solutions to system~\eqref{1.1}, which, to the best of our knowledge, only includes the qualitative analysis performed in \cite{HN25} for constant choices of $f$, as well as periodicity properties for the associated ODE system \eqref{1.3} in the case $a=0$ and a time-periodic $\tilde{f}(t)$. In particular, \cite{HN25} proves that the balance between the logistic growth rate $r$ and the external supply of the chemical $f$ governs the linearized dynamics of the system.
\\\\
More specifically, for a constant source, $f(x,t) \equiv f >0$, if $f > r$, it is proven that the only nonnegative spatially homogeneous steady state $(0,f)$ is locally asymptotically stable, while if $f < r$, there exists a secondary nonnegative equilibrium $(u_*, v_*) := \left(\displaystyle\frac{r-f}{r+a}, \frac{r(f+a)}{r+a} \right)$. In this case, $(0,f)$ is unstable, while $(u_*, v_*)$ is locally asymptotically stable. The local dynamics were at times preserved globally, while certain parameter gave rise to aggregation phenomena in the form of spike patterns. 
\\\\
For the case in which $a=0$ and $\tilde{f}(t)$ is a time-periodic function of period $T>0$, it is shown that if $r > \bar{\tilde{f}} := \displaystyle \frac{1}{T} \int_0^T \tilde{f}(t) \hspace{0.1 cm} dt$, then system \eqref{1.3} admits a periodic solution of period $T$, thus inheriting the periodicity of the source $\tilde{f}$. The convergence of solutions $(u,v)$ to $(\tilde{u},\tilde{v})$ was only tested numerically with positive results.
\\\\
In this way, the present work provides an analytical proof of the convergence of solutions to the full PDE system \eqref{1.1} to its associated ODE system \eqref{1.3}, that holds as well for the more interesting case $a>0$, when the second equation is not decoupled.
\\\\
In particular, for the proof it is assumed that $f = f(x,t)$ and $\tilde{f} = \tilde{f}(t)$ are nonnegative functions satisfying the following properties
\begin{equation} \label{1.2}
\begin{split}
    f \in L^\infty \big(\Omega \times(0, \infty) \big), 
        \\\\
         \int_0^\infty ||f-\tilde{f}||_{L^1(\Omega)} \hspace{0.1 cm} dt < \infty,
\end{split}
\end{equation}
which are also considered in \eqref{ant-restricciones} for system \eqref{ant}. However, for our proof, we need to also assume a stronger convergence of $f$ to $\tilde{f}$ in the sense that
\begin{equation}\label{1.2.1}
\int_0^\infty ||f - \tilde{f}||^2_{L^2(\Omega)} \hspace{0.1 cm} dt + \int_0^\infty \Big|\Big| \int_\Omega f - \tilde{f}\Big|\Big|^2_{L^2(\Omega)}  \hspace{0.1 cm} dt < \infty.
\end{equation}
In contrast, the second hypothesis for $f$ in \eqref{ant-restricciones} regarding the integrability in time of $||\nabla f||_{L^2(\Omega)}^2$ is not needed for our study. 
\\\\
In the most natural case in which $\tilde{f}$ is periodic, $f$ represents a source term that is asymptotically periodic in time, globally bounded, and whose spatial heterogeneities decay in time, in the sense that the second part of \eqref{1.2} and \eqref{1.2.1} hold. This models a persistent global periodic source with certain localized spatial perturbations that are however controlled in the long term. 
\\\\
For instance, $f$ can be taken as
$$f(x,t) = \tilde{f}(t) + p(t)q(x),$$
where $\tilde{f}\in L^\infty\big((0, \infty)\big)$ is time periodic, $q \in L^1(\Omega) \cap L^\infty(\Omega)$ encapsulates the spatial heterogeneities, and $p \in L^1\big((0,\infty)\big) \cap L^\infty\big((0,\infty)\big)$ represents the time decay of such heterogeneities.
\\\\
In nature, there are various examples of chemotactic processes exhibiting such periodic dynamics. For instance, in the presence of the chemoattractant cyclic AMP, in \cite{dict} the movement of \textit{Dictyostelium discoideum} amoebae was analyzed, concluding that the average velocity towards their center of aggregation is periodic in time. An example more closely related to our scenario was reported in \cite{Freq}, in which \textit{E. coli} cells were exposed to spatially and temporally varying attractant sources. It was shown that at low source frequencies, the bacterial population density oscillates in synchrony with the periodically supplied attractant. 
\\\\
These observations provide a biological motivation for analyzing the long-time behavior of Keller-Segel-type systems subject to periodic external inputs, and for identifying the conditions under which the macroscopic variables inherit the periodicity of the source.
\\\\
The article is structured as follows. After this introduction, Section \ref{s2} covers the global existence of bounded solutions to the system, while Section \ref{s3} includes further estimates that allow us to prove that the solution $u$ converges to $\displaystyle \bar{u}(t) := \frac{1}{|\Omega|} \int_\Omega u(x,t) \hspace{0.1 cm} dx$. The main result of the article, Theorem \ref{t-per} is included in Section \ref{s4}, proving that indeed $\bar{u}(t)$ converges to $\tilde{u}(t)$, and proceeding similarly for $v$. A final discussion concerning the case in which $\tilde{f}$ is time periodic is included in Section \ref{s5}, with some comments regarding the conditions under which the solutions to the ODE system \eqref{1.3} inherit the periodicity of $\tilde{f}$.

\section{Existence and uniqueness of solutions}\label{s2}
In this section we present Theorem \ref{t-exst}, the main result concerning the global existence of bounded solutions to system \eqref{1.1}. For brevity reasons, as the proof is similar to others and mainly follows the steps in \cite{Winkler10}, we only provide an outline of it. 
\begin{theorem}\label{t-exst}
    Let $\Omega \subset \mathbb{R}^n$, $n \geq 1$ be a convex, bounded domain with smooth boundary and assume that $f \in L^\infty \big(\Omega \times (0, \infty)\big)$, $u_0 \in C^0(\bar{\Omega})$ and $v_0 \in W^{1,\infty}(\Omega)$. Then, for all $\chi \in \mathbb{R}$, $D>0$ and $a\geq 0$ there exists $r_0 >0$ such that if $r\geq r_0$, problem \eqref{1.1} with boundary and initial conditions \eqref{1.1.2} admits a uniquely determined nonnegative global classical solution $(u,v)$ with
    $$
    ||u||_{L^\infty(\Omega)} + ||v||_{W^{1,\infty}(\Omega)} < c, \quad \text{for all } t>0,
    $$
    for a certain $c>0$ independent of $t$.
\end{theorem}
\textit{Sketch of the proof} As above mentioned, the proof can be carried out similarly as in \cite{Winkler10}. Firstly, the local existence of solutions in an interval $(0,T_{\max})$ is standard following the results by Amann \cite{amann}. Then, after appropriate estimates for the quantity $\displaystyle \sum_{k=0}^m b_k \int_\Omega u^k |\nabla v|^{2m-2k}$, for conveniently selected weights $b_k$, time-independent bounds for $u$ in $L^m(\Omega)$ and for $\nabla v$ in $L^{2m}(\Omega)$ can be obtained for arbitrarily large $m \in \mathbb{N}$. The main difference with respect to \cite{Winkler10} lies on the conditions for $r_0$ based on previois estimates that depend on $||f||_{L^\infty(\Omega \times (0,\infty))}$. Lastly, standard semigroup theory allows to obtain the desired bounds for $u$ and $v$.

\section{Further estimates}\label{s3}
Next, we derive some estimates that play a crucial role when analyzing the asymptotic behavior of the solutions. The first result entails a lower bound for the total mass of $u$.
\begin{lemma}\label{l3.1}
Let $f$ be such that $||f||_{L^\infty(\Omega \times (0,\infty))} < \infty$. Then, if $r>0$ is large enough, there exists $c>0$ satisfying 
$$\int_\Omega u \geq c, \quad \text{for all } t>0.$$
\end{lemma}
\begin{proof}
Firstly, the assumptions on $f$ and $r$ ensure that by Theorem \ref{t-exst} there exists a unique global solution $(u,v)$ to system \eqref{1.1}, satisfying  $u \in L^\infty\big(\Omega \times (0,\infty)\big)$ and 
$v \in L^\infty \big((0,\infty); ~W^{1,\infty}(\Omega) \big)$. 
\\\\
The main part of the proof follows the method introduced by Mizukami and Yokota in \cite{mz-yk}, Lemma 4.2. To this end, we take an arbitrary $\beta \in (0,1)$ satisfying
\begin{equation} \label{3.01}
     0<\beta < \min \left \{ 1,  \frac{4 \left( r - ||v||_{L^\infty\left(\Omega \times (0,\infty)\right)} \right)}{D \chi^2 ||\nabla v||^2_{L^\infty(\Omega \times (0,\infty))}}  - 1 \right \},
\end{equation}
and consider $\displaystyle \int_\Omega u^{-\beta} $. In particular, the choice of $r>0$ must be large enough to satisfy
\begin{equation}\label{3.0}
r> \frac{D \cdot  \chi^2 \cdot ||\nabla v||^2_{L^\infty\left(\Omega \times (0,\infty)\right)}}{4}+ ||v||_{L^\infty\left(\Omega \times (0,\infty)\right)},
\end{equation}
ensuring that indeed $\beta >0$. Integrating by parts we obtain
\begin{equation}\label{3.1}
\begin{split}
    \frac{d}{dt} \int_\Omega u^{-\beta} =& - \beta \int_\Omega u^{-\beta - 1} \Bigg ( D\Delta u - \chi \nabla \cdot(u \nabla v) + ru \left(1-u - \frac{v}{r} \right) \Bigg) \\
    =& -D \beta(\beta+1) \int_\Omega  u^{-\beta -2} |\nabla u|^2 - \chi \beta(\beta+1) \int_\Omega u^{-\beta - 1} \nabla u \cdot \nabla v \\
     &- r\beta \int_\Omega u^{-\beta} \left(1-u - \frac{v}{r} \right), \quad \text{ for all t>0}.
\end{split}
\end{equation}
The term involving the product of the gradients can be bounded using Young's inequality as follows
\begin{equation} \label{3.2}
\begin{split}
    - \chi \beta(\beta+1) \int_\Omega u^{-\beta - 1} \nabla u \cdot \nabla v  \leq \chi \beta(\beta+1) \int_\Omega u^{-\beta - 1} |\nabla u| |\nabla v| \\
    \leq D \beta (\beta +1)  \int_\Omega u^{-\beta -2} |\nabla u|^2 + \frac{D \beta (\beta+1) \chi^2}{4} \int_\Omega u^{-\beta} |\nabla v|^2,
\end{split}
\end{equation}
for all $t>0$. Notice that the first term appearing on the right hand side here cancels the corresponding one in \eqref{3.1}. We can moreover estimate the final term in \eqref{3.2} by making use of the $L^\infty \big(\Omega \times (0,\infty)\big)$ bound for $\nabla v$. In particular, for all $t>0$ we have
\begin{equation} \label{3.3}
    \frac{D \beta (\beta+1) \chi^2}{4} \int_\Omega u^{-\beta} |\nabla v|^2 \leq  
    ||\nabla v||^2_{L^\infty(\Omega \times (0,\infty)} \cdot 
    \frac{D \beta (\beta+1) \chi^2}{4} \int_\Omega u^{-\beta} .
\end{equation}
Next, concerning the last term appearing in \eqref{3.1}, we have
$$
-r \beta \int_\Omega u^{-\beta} \left(1-u - \frac{v}{r} \right) = -r \beta \int_\Omega u^{-\beta} + r \beta \int_\Omega u^{1-\beta} + \beta \int_\Omega u^{-\beta} v.
$$
Relying this time on the global boundedness of $u$, the second term on the right hand side of the above expression is such that
\begin{equation}\label{3.4}
    r \beta \int_\Omega u^{1-\beta} \leq r \beta |\Omega| ||u||_{L^\infty\left(\Omega \times (0,\infty)\right)}^{1-\beta}, \quad \text{for all } t>0,
\end{equation}
and for the last one we have
\begin{equation}\label{3.5}
    \beta \int_\Omega u^{-\beta} v \leq \beta ||v||_{L^\infty\left(\Omega \times (0,\infty)\right)} \int_\Omega u^{-\beta} \quad \text{for all } t>0.
\end{equation}
In this way, by combining estimates \eqref{3.2}-\eqref{3.5}, we can rewrite \eqref{3.1} directly as
\begin{equation} \label{3.6}
    \frac{d}{dt} \int_\Omega u^{-\beta} +  c_1 \int_\Omega u^{-\beta}\leq c_2, \quad \text{for all } t>0,
\end{equation}
where
\begin{equation*}
    \begin{split}
        c_1 &:= \beta \left(r -||v||_{L^\infty\left(\Omega \times (0,\infty)\right)} - ||\nabla v||^2_{L^\infty(\Omega \times (0,\infty)} \cdot \frac{ D (\beta+1) \chi^2}{4}\right), 
        \\\\
        c_2 &:=  r \beta |\Omega| ||u||_{L^\infty\left(\Omega \times (0,\infty)\right)}^{1-\beta} >0.
    \end{split}
\end{equation*}
The choice of $\beta$ in \eqref{3.01}, combined with assumption \eqref{3.0} ensures that indeed $c_1 >0$ as well. Thus, the differential inequality \eqref{3.6} implies that there exists $c_\beta >0$ such that
\begin{equation}\label{3.7}
    \int_\Omega u^{-\beta} \leq c_\beta, \quad \text{for all } t>0.
\end{equation}
Next, combining the the Cauchy Schwarz inequality with the obtained bound \eqref{3.7}, we have for all $t>0$
$$
|\Omega| = \int_\Omega \frac{u^{\beta/2}}{u^{\beta/2}} \leq \left[ \int_\Omega u^\beta \right]^{\frac{1}{2}} \cdot \left[ \int_\Omega u^{-\beta} \right]^{\frac{1}{2}} \leq c_\beta^{\frac{1}{2}} \cdot \left[ \int_\Omega u^\beta \right]^{\frac{1}{2}}.
$$
Lastly, as $\beta \in (0,1)$, using Hölder's inequality with exponents $p = \beta^{-1}>1$, $q = (1-\beta)^{-1} > 1$, it follows that
$$
|\Omega|^2 \cdot c_{\beta}^{-1} \leq \int_\Omega u^\beta  = \int_\Omega \left(u^\beta \cdot 1\right)\leq \left ( \int_\Omega u \right)^\beta |\Omega|^{1-\beta}, \quad \text{for all } t>0,
$$
which finishes the proof upon taking $c:= \left(|\Omega|^{1+\beta} \cdot c_{\beta}^{-1}\right)^{1/\beta}$.
\end{proof}
Next, we prove a similar result concerning a lower bound for the solution $\tilde{u}$ to the associated ODE system. In this case it is again a threshold for the logistic growth rate $r$ which drives the dynamics of the system. 
\begin{lemma}\label{l3.2}
    Let $\tilde{u}$ be the solution to system \eqref{1.3} with initial values \eqref{1.4} and assume that 
    \begin{equation}\label{3.8}
        r > \max \big \{ \tilde{v}(0), ~ a \cdot \max\{\tilde{u}(0), 1\} + ||\tilde{f}||_{L^\infty (0, \infty)} \big\}.
    \end{equation}
    Then, there exists $c>0$ such that 
    $$\tilde{u}(t) > c, \quad \text{for all } t >0.$$
\end{lemma}
\begin{proof}
From the first equation in system \eqref{1.3} we have
\begin{equation}\label{3.9}
    \frac{1}{r \tilde{u}} \cdot  \frac{~d \tilde{u}}{dt} + \tilde{u} = 1- \frac{\tilde{v}}{r}, \quad \text{for all } t>0.
\end{equation}
The main element in the proof is based on the fact that under hypothesis \eqref{3.8}, the solution $\tilde{v}$ is such that $\displaystyle 1 - \frac{\tilde{v}}{r} >0$. To prove this, we start by noting that $\tilde{u}$ satisfies
$$
\frac{~d \tilde{u}}{dt} = r\tilde{u}(1-\tilde{u}) - \tilde{u}\tilde{v} \leq r\tilde{u}(1-\tilde{u}),
$$
for all $t>0$, which directly implies that $\tilde{u}$ is bounded from above by the solution of the logistic equation, and in particular 
\begin{equation}\label{3.10}
    \tilde{u}(t) \leq \max\{ \tilde{u}(0), 1\}, \quad \text{for all } t >0.
\end{equation}
Substituting the bound into the second equation in \eqref{1.3} yields
$$
\frac{~d \tilde{v}}{dt} = a\tilde{u} - \tilde{v} + \tilde{f} \leq \left[ a  \cdot \max\{ \tilde{u}(0), 1\} + ||\tilde{f}||_{L^\infty (0, \infty)} \right] - \tilde{v},
$$
for all $t >0$. Solving the differential inequality, one obtains
\begin{equation}\label{3.11}
    \tilde{v}(t) \leq \max \big \{\tilde{v}(0),  a  \cdot \max\{ \tilde{u}(0), 1\} + ||\tilde{f}||_{L^\infty (0, \infty)} \big \}, \quad \text{for all } t >0.
\end{equation}
Thus, if $\displaystyle r > \max \big \{\tilde{v}(0),  a  \cdot \max\{ \tilde{u}(0), 1\} + ||\tilde{f}||_{L^\infty (0, \infty)} \big \}$, \eqref{3.9} implies that there exists $\varepsilon >0$ such that
\begin{equation} \label{3.12}
 \frac{1}{r \tilde{u}} \cdot  \frac{~d \tilde{u}}{dt} + \tilde{u}  > \varepsilon, \quad \text{for all } t >0.
\end{equation}
The remaining of the proof is standard. By defining $\displaystyle \tilde{\varepsilon} := \frac{1}{2} \min \{\tilde{u}(0), \varepsilon\} >0$, we prove that $\tilde{u}(t) > c$ for all $t>0$. 
\\\\
We proceed by contradiction by supposing there exists $t_0 > 0$ such that $\tilde{u}(t)>\tilde{\varepsilon}$ for all $t \in (0,t_0)$, while also satisfying $\tilde{u}(t_0) = \tilde{\varepsilon}$ and $\displaystyle \frac{~d\tilde{u}}{dt}(t_0) \leq 0$. In this way, at $t = t_0$ \eqref{3.12} turns into
$$
\frac{1}{r \tilde{\varepsilon}} \cdot \frac{~d\tilde{u}}{dt}(t_0) + \tilde{\varepsilon} > \varepsilon.
$$
Thus, as $\varepsilon > \tilde{\varepsilon}$ we have
$$
\frac{~d\tilde{u}}{dt}(t_0) > r \tilde{\varepsilon} ( \varepsilon - \tilde{\varepsilon}) >0,
$$
which contradicts our initial assumption that $\displaystyle \frac{~d\tilde{u}}{dt}(t_0) \leq 0$. Thus, the proof is complete by taking $c := \tilde{\varepsilon}$.
\end{proof}
Next, we introduce the nonnegative function
\begin{equation}\label{3.13}
    k_1(t) := \int_\Omega \left(u(x,t) - \int_\Omega u(x,t)~ dx  \right)^2dx,
\end{equation}
which verifies the integrability property proved in the following lemma.
\begin{lemma}\label{l3.3}
    Let $f$ be such that $||f||_{L^\infty(\Omega \times (0,\infty))} < \infty$, and that 
    \begin{equation}\label{hipk3}
        \begin{split}
       & \int_0^\infty ||f - \tilde{f}||_{L^1(\Omega)} \hspace{0.1 cm} dt < \infty \\
         &   \int_0^\infty ||f - \tilde{f}||^2_{L^2(\Omega)} \hspace{0.1 cm} dt + \int_0^\infty \Big|\Big| \int_\Omega f - \tilde{f}\Big|\Big|^2_{L^2(\Omega)}  \hspace{0.1 cm} dt < \infty.
        \end{split}
    \end{equation}
    for some spatially homogeneous function $\tilde{f}$. Then, if $r>0$ is large enough, the function $k_1$ defined in \eqref{3.13} verifies
    $$
    \int_0^\infty k_1(t) ~dt \leq c < \infty,
    $$
    for a certain $c>0$.
\end{lemma}
\begin{proof}
First, the assumption on $r$ being large enough grants the applicability of Theorem \ref{t-exst} and Lemmas \ref{l3.1} and \ref{l3.2}. To prove the result, we first decompose $v$ into two components, $w$ and $z$, each one accounting for one of the terms in $\displaystyle \frac{\partial v}{\partial t}$. 
\\\\
More precisely, we let $w$ and $z$ be the solutions of
\begin{equation}\label{3.14}
    \begin{cases}
        \displaystyle \frac{\partial w}{\partial t} = \Delta w + au - w, \quad & x \in \Omega, ~t > 0,  \\\\
        \displaystyle \frac{\partial z}{\partial t} = \Delta z + f - z, \quad & x \in \Omega, ~t > 0, 
    \end{cases}
\end{equation}
with initial values $w(x,0) = z(x,0) = \displaystyle \frac{1}{2} v(x,0)$, and Neumann homogeneous boundary conditions, clearly having
$$v = w + z.$$
In this way, the dynamics of $w$ capture the production term $+au$, while the source term $+f$ is handled solely by the equation for $z$. Naturally, the existence theory developed for $v$ in Theorem \eqref{t-exst} extends to both $w$ and $z$, resulting in their global boundedness.
\\\\
The equation for $u$ is now given by
\begin{equation}\label{3.15}
     \frac{\partial u}{\partial t} = D\Delta u + \chi \nabla \cdot \Big(u \nabla (w+z) \Big) + ru\left(1- u - \frac{w}{r} - \frac{z}{r}\right),
\end{equation}
for $ x \in \Omega, ~t > 0$. Moreover, the associate ODE system \eqref{1.3} now becomes
\begin{equation}\label{3.16}
\begin{cases}
      \displaystyle  \frac{d \tilde{u}}{dt~} = r \tilde{u} \left(1 - \tilde{u} - \frac{\tilde{w}}{r} - \frac{\tilde{z}}{r} \right), \quad & t>0, \\\\
      \displaystyle  \frac{d \tilde{w}}{dt~} = a\tilde{u} - \tilde{w}, \quad & t>0,\\\\
      \displaystyle  \frac{d \tilde{z}}{dt~} = \tilde{f} - \tilde{z}, \quad & t>0,
\end{cases}
\end{equation}
with initial values $\displaystyle \tilde{u}(0) = \int_\Omega u_0(x) \hspace{0.1 cm} dx$, $\displaystyle \tilde{w}(0) = \frac{1}{2}\int_\Omega v_0(x) \hspace{0.1 cm} dx$, and $\displaystyle \tilde{z}(0) = \frac{1}{2}\int_\Omega v_0(x) \hspace{0.1 cm} dx$. 
\\\\
We start by defining the following functionals
\begin{equation}\label{3.24}
\begin{split}
    &F_1(t) := \int_\Omega \frac{u(x,t)}{\tilde{u}(t)} ~dx - 1 + \ln \tilde{u}(t) - \int_\Omega \ln u(x,t) ~dx, 
    \\\\
    &F_2(t) := \ln \left(\int_\Omega u(x,t) ~dx \right) - \ln \tilde{u}(t),
\end{split}
\end{equation}
with the same structure as others used in the study of chemotaxis systems, such as in \cite{bw1, per1}. Given that $|\Omega|=1$, $F_1$ verifies
$$
F_1(t) = \int_\Omega h\left(\frac{u}{~\tilde{u}} \right), \quad \text{for } h(s):= s-1-\ln s,
$$
with $h(s)\geq 0$ for all $s>0$, having $\lim_{s \to 0^+} h(s) = + \infty$. Thus, $F_1(t) \geq 0$. With respect to $F_2$, Lemma \ref{l3.1} and the fact that $\tilde{u}(t) \leq \max\{ \tilde{u}(0), 1\}$ for all $t>0$ grant the existence of a $c_3>0$ such that $F_2(t) \geq -c_3 > - \infty$ for all $t>0$.
\\\\
We first compute the time derivative of $F_2$. We start by integrating equation \eqref{3.15} over $\Omega$, yielding
\begin{equation}\label{3.17}
    \begin{split}
      \frac{1}{r} \frac{d}{dt} &\int_\Omega u = \frac{1}{r} \left[D\int_\Omega \Delta u + \chi \int_\Omega \nabla \cdot \Big(u \nabla (w+z) \Big) \right] + \int_\Omega u\left(1- u - \frac{z}{r}\right) -  \frac{1}{r}\int_\Omega u w
        \\\\
        & = \int_\Omega \left(u - \int_\Omega u + \int_\Omega u \right) \left(1- u - \frac{z}{r}\right) - \frac{1}{r}\int_\Omega u w 
        \\\\
        & =  \int_\Omega \left(u - \int_\Omega u \right) \left(1- u - \frac{z}{r}\right) + \int_\Omega \left(\int_\Omega u \right)\left(1- u - \frac{z}{r}\right) - \frac{1}{r}\int_\Omega u w
        \\\\
        &  = \int_\Omega \left[ \left(u - \int_\Omega u \right) \left(1- u + \int_\Omega u - \int_\Omega u - \frac{z}{r} +  \frac{\tilde{z}}{r} -  \frac{\tilde{z}}{r} \right) \right]
        \\\\
        & + \left(\int_\Omega u \right) \cdot \int_\Omega \left(1- u - \frac{z}{r}\right)  - \frac{1}{r}\int_\Omega u w 
        \\\\
        & = \int_\Omega \left[ \left(u - \int_\Omega u \right)  \left(\int_\Omega u - u \right) \right] + \int_\Omega \left[ \left(u - \int_\Omega u \right) \left(  \frac{\tilde{z}}{r} -  \frac{z}{r} \right) \right]
        \\\\
        & + \int_\Omega \left[ \left(u - \int_\Omega u \right) \left(1 - \frac{\tilde{z}}{r} - \int_\Omega u \right) \right] +  \left(\int_\Omega u \right) \cdot \int_\Omega \left(1- u - \frac{z}{r}\right)  - \frac{1}{r}\int_\Omega u w,
    \end{split}
\end{equation}
for all $t>0$. Simplifying the terms in the second-to-last line of the above expression, we obtain for all $t >0$
\begin{equation}\label{3.18}
    \begin{split}
        \int_\Omega & \left[ \left(u - \int_\Omega u \right) \left(1 - \frac{\tilde{z}}{r} - \int_\Omega u \right) \right] +  \left(\int_\Omega u \right) \cdot \int_\Omega \left(1- u - \frac{z}{r}\right)  
        \\\\
        & = \int_\Omega u \left(1 - \frac{\tilde{z}}{r} - \int_\Omega u \right)  - \left(\int_\Omega u \right) \cdot \int_\Omega \left(1 - \frac{\tilde{z}}{r} - \int_\Omega u \right)
        \\\\
        &  +  \left(\int_\Omega u \right) \cdot \int_\Omega \left(1- u - \frac{z}{r}\right) 
        \\\\
        & = \int_\Omega u \left(1 - \frac{\tilde{z}}{r} - \int_\Omega u \right) + \left( \int_\Omega u \right) \cdot \int_\Omega \left(\frac{\tilde{z}}{r} - \frac{z}{r} \right),
    \end{split}
\end{equation}
where we recall that we assumed that $|\Omega| = 1$. Hence, by combining \eqref{3.17} and \eqref{3.18} we arrive at
\begin{equation}\label{3.19}
\begin{split}
     \frac{1}{r} \frac{d}{dt} &\int_\Omega u = \int_\Omega \left[ \left(u - \int_\Omega u \right)  \left(\int_\Omega u - u \right) \right] + \int_\Omega \left[ \left(u - \int_\Omega u \right) \left(  \frac{\tilde{z}}{r} -  \frac{z}{r} \right) \right]
     \\\\
     & \int_\Omega u \left(1 - \frac{\tilde{z}}{r} - \int_\Omega u \right) + \left( \int_\Omega u \right) \cdot \int_\Omega \left(\frac{\tilde{z}}{r} - \frac{z}{r} \right) - \frac{1}{r} \int_\Omega uw
     \\\\
     & = - k_1(t) +  \int_\Omega \left[ \left(u - \int_\Omega u \right) \left(  \frac{\tilde{z}}{r} -  \frac{z}{r} \right) \right] +  \int_\Omega u \left(1 - \frac{\tilde{z}}{r} - \int_\Omega u \right) 
     \\\\
     & + \Bigg|\Bigg| \frac{\tilde{z}}{r} - \frac{z}{r} \Bigg|\Bigg|_{L^1(\Omega)} \cdot \left( \int_\Omega u \right) -  \frac{1}{r} \int_\Omega uw, \quad \text{for all } t>0.
\end{split}
\end{equation}
For the second term appearing on the right hand side of \eqref{3.19}, we consider $\delta >0$ given by 
\begin{equation}\label{delta}
    \delta = \frac{1}{2} \min \left\{1, 1- \frac{\varepsilon_u}{r} \left(\frac{\chi^2 a^2}{8D} - \frac{a}{||u||_{L^\infty(\Omega \times(0,\infty))}} \right) \right \},
\end{equation}
where $\varepsilon_u >0$ is the lower bound for $\displaystyle \int_\Omega u$ provided by Lemma \ref{l3.1}. Again the assumption on $r$ being large enough allows us to grant that $\delta>0$.
\\\\
With the value of $\delta$, using Young's inequality we have
\begin{equation}\label{3.20}
    \begin{split}
        \int_\Omega & \left[ \left(u - \int_\Omega u \right) \left(  \frac{\tilde{z}}{r} -  \frac{z}{r} \right) \right]  \leq \delta \int_\Omega \left(u - \int_\Omega u \right)^2 + \frac{1}{4 \delta} \int_\Omega  \left(  \frac{\tilde{z}}{r} -  \frac{z}{r} \right)^2
        \\
        & \leq \delta k_1 + \frac{||\tilde{z}-z||_{L^\infty(\Omega)}}{4\delta} \cdot  \Bigg|\Bigg| \frac{\tilde{z}}{r} - \frac{z}{r} \Bigg|\Bigg|_{L^1(\Omega)} =: ~\delta k_1 + c_1(\delta) \Bigg|\Bigg| \frac{\tilde{z}}{r} - \frac{z}{r} \Bigg|\Bigg|_{L^1(\Omega)},
    \end{split}
\end{equation}
for all $t>0$, where we used the global boundedness of $z$ and $\tilde{z}$, as a direct consequence of the boundedness of $f$ and $\tilde{f}$. Therefore, substituting in \eqref{3.19} we obtain
\begin{equation}\label{3.21}
\begin{split}
    \frac{1}{r} \frac{d}{dt}  \int_\Omega u & \leq  -(1-\delta) k_1 + \int_\Omega u \left(1 - \frac{\tilde{z}}{r} - \int_\Omega u \right)  
    \\ &+ c_2(\delta) \Bigg|\Bigg| \frac{\tilde{z}}{r} - \frac{z}{r} \Bigg|\Bigg|_{L^1(\Omega)} - \frac{1}{r} \int_\Omega uw, \quad \text{for all }t>0,
\end{split}
\end{equation}
where $c_2(\delta) := c_1(\delta) + ||u||_{L^\infty(\Omega)}>0$. Next, dividing this last inequality by $\displaystyle \int_\Omega u$, we have for all $t>0$
\begin{equation}\label{3.22}
    \begin{split}
    \frac{d}{dt}   \ln  \left(\int_\Omega u \right) \leq r \Bigg\{ \frac{-(1-\delta)k_1}{\displaystyle \int_\Omega u} + \frac{\displaystyle\int_\Omega u \left(1 - \frac{\tilde{z}}{r} - \int_\Omega u \right) }{\displaystyle \int_\Omega u} +  \frac{c_2(\delta)}{\displaystyle \int_\Omega u} \cdot \Bigg|\Bigg| \frac{\tilde{z}}{r} - \frac{z}{r} \Bigg|\Bigg|_{L^1(\Omega)}~ \Bigg \}
    \\
      - \frac{ \displaystyle \int_\Omega uw}{\displaystyle \int_\Omega u} = r \Bigg\{ \frac{-(1-\delta)k_1}{\displaystyle \int_\Omega u} + 1 + \frac{\tilde{z}}{r} - \int_\Omega u \Bigg \}  +  \frac{c_2(\delta)}{\displaystyle \int_\Omega u} \cdot || \tilde{z} -z  ||_{L^1(\Omega)} - \frac{ \displaystyle \int_\Omega uw}{\displaystyle \int_\Omega u}.
    \end{split}
\end{equation}
In turn, $\tilde{u}$ verifies
$$
    \frac{d}{dt} \ln \tilde{u} = r \left( 1 - \tilde{u} - \frac{\tilde{w}}{r} - \frac{\tilde{z}}{r}\right), \quad \text{for all } t >0,
$$
and thus
\begin{equation}\label{3.23}
    \begin{split}
   \frac{d}{dt} F_2 &= \frac{d}{dt} \left( \ln  \left(\int_\Omega u \right)  - \ln \tilde{u} \right) \leq r \left( - \frac{(1-\delta) k_1}{\displaystyle \int_\Omega u} + \frac{\tilde{w}}{r} + \tilde{u} - \int_\Omega u \right) \\
    &+ \frac{ c_2(\delta)}{\displaystyle \int_\Omega u} || \tilde{z} -z  ||_{L^1(\Omega)} - \frac{ \displaystyle \int_\Omega u w}{\displaystyle\int_\Omega w}, \quad \text{for all } t >0.
    \end{split}
\end{equation}
With respect to $F_1$, by means of an integration by parts followed by Young's inequality, we have
\begin{equation}\label{3.25}
\begin{split}
    \frac{d}{dt} F_1 &= \frac{d}{dt} \left( \int_\Omega \frac{u}{\tilde{u}}\right) +\frac{~d \tilde{u}}{dt} \cdot \frac{1}{~\tilde{u}} - \int_\Omega \frac{\partial u}{\partial t} \cdot \frac{1}{u} 
    \\\\
    & = \frac{d}{dt} \left( \int_\Omega \frac{u}{\tilde{u}}\right)  + r \left(1-\tilde{u}-\frac{\tilde{w}}{r} - \frac{\tilde{z}}{r} \right) - \int_\Omega \frac{D |\nabla u|^2}{u^2}
    \\\\
    & -  \int_\Omega \chi \frac{\nabla u \cdot \nabla w}{u} - \int_\Omega \chi \frac{\nabla u \cdot \nabla z}{u} + r \int_\Omega u \left(1-u- \frac{w}{r}-\frac{z}{r} \right) 
    \\\\
    & \leq \frac{d}{dt} \left( \int_\Omega \frac{u}{\tilde{u}}\right) - D \int_\Omega \frac{ |\nabla u|^2}{u^2} + \frac{D}{2} \int_\Omega \frac{ |\nabla u|^2}{u^2} + \frac{\chi^2}{2D}\int_\Omega |\nabla w|^2
    \\\\
    & + \frac{D}{2}\int_\Omega \frac{ |\nabla u|^2}{u^2}  + \frac{\chi^2}{2D} \int_\Omega |\nabla z|^2 + r \int_\Omega u  - \tilde{u} + \int_\Omega w - \tilde{w} + \int_\Omega z - \tilde{z}
    \\\\
    & = \frac{d}{dt} \left( \int_\Omega \frac{u}{\tilde{u}}\right) +  r \int_\Omega u  - \tilde{u} + \int_\Omega w - \tilde{w} + \int_\Omega z - \tilde{z} 
    \\\\
    & + \frac{\chi^2}{2} \left[\int_\Omega |\nabla w|^2 + \int_\Omega |\nabla z|^2  \right], \quad \text{for all } t >0.
    \end{split}
\end{equation}
The next step is to bound the last term on the right hand side involving the gradients of $w$ and $z$. For the first one, multiplying the first equation in \eqref{3.14} by $w$ and integrating by parts over $\Omega$, we obtain
$$
\frac{1}{2} \frac{d}{dt} \int_\Omega w^2 = - \int_\Omega |\nabla w|^2 + a \int_\Omega u w - \int_\Omega w^2, \quad \text{for all } t >0.
$$
As $|\Omega| = 1$, for every function $g\in L^2(\Omega)$, the well known variance identity holds
\begin{equation}\label{3.id1}    
\int_\Omega \left(g - \int_\Omega g \right)^2 = \int_\Omega g^2 - \left( \int_\Omega g\right)^2,
\end{equation}
and thus we arrive at
\begin{equation}\label{3.26}
    \begin{split}
        \int_\Omega |\nabla w|^2 &=  a \int_\Omega u w - \int_\Omega w^2- \frac{1}{2} \frac{d}{dt} \int_\Omega w^2 
        \\
        &= a \int_\Omega u w - \int_\Omega w^2 - \frac{1}{2} \frac{d}{dt} \int_\Omega \left( w - \int_\Omega w \right)^2 - \int_\Omega w \cdot \int_\Omega \frac{\partial w}{\partial t}
        \\
        & = a \int_\Omega u w - \int_\Omega w^2 - \frac{1}{2} \frac{d}{dt} \int_\Omega \left( w - \int_\Omega w \right)^2
        \\
        & - \int_\Omega w \cdot \int_\Omega \big( \Delta w + au - w\big)
        \\
        & = a \int_\Omega uw - a \int_\Omega u \cdot \int_\Omega w  - \frac{1}{2} \frac{d}{dt} \int_\Omega \left( w - \int_\Omega w \right)^2 
        \\
        & + \left(\int_\Omega w \right)^2 - \int_\Omega w^2, \quad \text{for all } t>0.
    \end{split}
\end{equation}
Next, one can check that the following identity holds for all $t>0$
$$
a \int_\Omega uw - \int_\Omega w^2 = \frac{a^2}{4}\left(\int_\Omega u \right)^2 - \int_\Omega \left(w-\frac{a}{2}u \right )^2 - \frac{a^2}{4} \left[ \left(\int_\Omega u \right)^2 - \int_\Omega u^2 \right], 
$$
so after substituting into \eqref{3.26} and completing the square with the remaining terms, we obtain
\begin{equation}\label{3.27}
    \begin{split}
        \int_\Omega & |\nabla w|^2 =  - \frac{1}{2} \frac{d}{dt} \int_\Omega \left( w - \int_\Omega w \right)^2  - \int_\Omega \left(w-\frac{a}{2}u \right )^2 - \frac{a^2}{4} \left[ \left(\int_\Omega u \right)^2 - \int_\Omega u^2 \right]
        \\
        & - a \int_\Omega u \cdot \int_\Omega w + \left(\int_\Omega w \right)^2 + \frac{a^2}{4}\left(\int_\Omega u \right)^2
        \\
        & = - \frac{1}{2} \frac{d}{dt} \int_\Omega \left( w - \int_\Omega w \right)^2  - \int_\Omega \left(w-\frac{a}{2}u \right )^2 - \frac{a^2}{4} \left[ \left(\int_\Omega u \right)^2 - \int_\Omega u^2 \right]
        \\
        & + \left[ \int_\Omega w - \frac{a}{2} \int_\Omega u \right]^2  =  - \frac{1}{2}\frac{d}{dt} \int_\Omega  \left( w - \int_\Omega w \right)^2  - \int_\Omega \left[ \left(w-\frac{a}{2}u \right )^2 - \int_\Omega \left( w-\frac{a}{2}u \right)\right]^2
        \\
        & - \frac{a^2}{4} \left[  \left(\int_\Omega u \right)^2 - \int_\Omega u^2 \right] \leq - \frac{1}{2}\frac{d}{dt} \int_\Omega \left( w - \int_\Omega w \right)^2+ \frac{a^2}{4} \int_\Omega \left(u - \int_\Omega u\right)^2 ,
    \end{split}
\end{equation}
for all $t>0$, where we relied again on identity \eqref{3.id1}. Recall that the final term, $\displaystyle \frac{a^2}{4} \int_\Omega \left(u - \int_\Omega u\right)^2$ is precisely $\displaystyle \frac{a^2}{4}  k_1$.
\\\\
Proceeding in a similar way, one can prove that
\begin{equation}\label{3.28}
    \int_\Omega |\nabla z|^2 \leq - \frac{1}{2} \frac{d}{dt} \int_\Omega \left(z - \int_\Omega z\right)^2 + \frac{1}{4} \int_\Omega \left(f - \int_\Omega f\right)^2 , \quad \text{for all } t >0.
\end{equation}
Combining bounds \eqref{3.27} and \eqref{3.28}, the time derivative of $F_1$ given in \eqref{3.25} can be rewritten as
\begin{equation}\label{3.29}
\begin{split}
    \frac{d}{dt} F_1 &=  \frac{d}{dt} \left( \int_\Omega \frac{u}{\tilde{u}}\right) +  r \int_\Omega u  - \tilde{u} + \int_\Omega w - \tilde{w} + \int_\Omega z - \tilde{z}  
    \\ 
    & + \frac{\chi^2}{2D} \left[ - \frac{1}{2} \frac{d}{dt} \int_\Omega \left( w - \int_\Omega w \right)^2 + \frac{a^2}{4} k_1  \right]
    \\
    & + \frac{\chi^2}{2D} \left[ - \frac{1}{2} \frac{d}{dt} \int_\Omega \left( z - \int_\Omega z \right)^2 + \frac{1}{4} \int_\Omega \left(f - \int_\Omega f\right)^2  \right], \quad \text{for all } t>0. 
\end{split}
\end{equation}
Thus, the combination of \eqref{3.29} with the bound for $\displaystyle \frac{d}{dt} F_2$ obtained in \eqref{3.23} yields, after some favorable cancellations
\begin{equation}\label{3.30}
    \begin{split}
    \frac{d}{dt} &(F_1+F_2)  \leq \frac{d}{dt} \left( \int_\Omega \frac{u}{\tilde{u}}\right) + \left( \frac{c_2(\delta)}{\displaystyle \int_\Omega u} + 1 \right) ||z-\tilde{z}||_{L^1(\Omega)} 
    \\&-  \left(\frac{r(1-\delta)}{\displaystyle \int_\Omega u} - \frac{\chi^2 a^2}{8D}\right) k_1+ \frac{\chi^2}{8D} \int_\Omega \left(f - \int_\Omega f\right)^2 + \int_\Omega w
    \\
    & - \frac{\chi^2}{4D} \frac{d}{dt} \left[\int_\Omega \left( w - \int_\Omega w \right)^2 + \int_\Omega \left( z - \int_\Omega z \right)^2 \right]
    - \frac{\displaystyle \int_\Omega uw}{\displaystyle \int_\Omega u},
    \end{split}
\end{equation}
for all $t>0$. Next, recall that $\varepsilon_u >0$ denotes the lower bound for $\displaystyle \int_\Omega u$ provided by Lemma \ref{l3.1}, and moreover, by Theorem \ref{t-exst} $\displaystyle \int_\Omega u(x,t) \hspace{0.1 cm} dx \leq ||u||_{L^\infty(\Omega \times (0,\infty))}<\infty$ for all $t>0$. as $|\Omega|=1$. Hence, we have
\begin{equation}\label{3.31}
    \begin{split}
        \int_\Omega w - & \frac{\displaystyle \int_\Omega uw}{\displaystyle \int_\Omega u}  =  \frac{\displaystyle \left(-\frac{d}{dt} \int_\Omega w + a\int_\Omega u\right) \cdot \left(\int_\Omega u\right) - \int_\Omega uw}{\displaystyle \int_\Omega u} 
        \\&\leq - \frac{d}{dt} \int_\Omega w + \frac{1}{\displaystyle \int_\Omega u}\cdot\left[ a\left(\int_\Omega u^2 - k_1\right) - \int_\Omega uw \right]
       \\& \leq - \frac{d}{dt} \int_\Omega w - \frac{a}{||u||_{L^\infty(\Omega \times (0,\infty))}} \cdot k_1 + \frac{||u||_{L^\infty(\Omega \times (0,\infty))}}{\varepsilon_u} \cdot \left[\int_\Omega \left(a u - w\right)\right]
       \\
       & = \left[\frac{||u||_{L^\infty(\Omega \times (0,\infty))}}{\varepsilon_u} - 1\right] \cdot \left(\frac{d}{dt}\int_\Omega w\right)   - \frac{a}{||u||_{L^\infty(\Omega \times (0,\infty))}} \cdot k_1,
    \end{split}
\end{equation}
for all $t>0$, where we relied once again in identity \eqref{3.id1} in the fact that
$$
\frac{d}{dt}\int_\Omega w = a \int_\Omega u - \int_\Omega w, \quad \text{for all } t>0.
$$
Direct substitution of \eqref{3.31} into \eqref{3.30} yields
\begin{equation*}
\begin{split}
\frac{d}{dt}&(F_1+F_2) +  \left(\frac{r(1-\delta)}{\displaystyle \int_\Omega u} - \frac{\chi^2 a^2}{8D} +  \frac{a}{||u||_{L^\infty(\Omega \times (0,\infty))}}\right) k_1 
\\
 &\leq \frac{d}{dt} \left( \int_\Omega \frac{u}{\tilde{u}}\right) + \left( \frac{c_2(\delta)}{\displaystyle \int_\Omega u} + 1 \right) ||z-\tilde{z}||_{L^1(\Omega)} 
+ \frac{\chi^2}{8D} \int_\Omega \left(f - \int_\Omega f\right)^2 \\
& - \frac{\chi^2}{4D} \frac{d}{dt} \left[\int_\Omega \left( w - \int_\Omega w \right)^2 + \int_\Omega \left( z - \int_\Omega z \right)^2 \right] + \left[\frac{||u||_{L^\infty(\Omega \times (0,\infty))}}{\varepsilon_u} - 1\right] \cdot \left(\frac{d}{dt}\int_\Omega w\right), 
\end{split}
\end{equation*}
for all $t>0$. Thus, a time integration over the interval $(0,\rho)$ for arbitrary $\rho>0$ yields
\begin{equation}\label{3.32}
\begin{split}
    A \int_0^\rho k_1(t) ~dt \leq & \Bigg[\int_\Omega \frac{u}{\tilde{u}} -F_1 -F_2 - \frac{\chi^2}{4D} \left(\int_\Omega \left( w - \int_\Omega w \right)^2 + \int_\Omega \left( z - \int_\Omega z \right)^2 \right) 
    \\
    &  + \left[\frac{||u||_{L^\infty(\Omega \times (0,\infty))}}{\varepsilon_u} - 1\right] \cdot \left(\frac{d}{dt}\int_\Omega w\right)\Bigg]_0^\rho + \frac{\chi^2}{8D} \int_0^\rho \int_\Omega \left(f - \int_\Omega f\right)^2
    \\
    & +\left( \frac{c_2(\delta)}{\displaystyle \int_\Omega u} + 1 \right) \int_0^\rho ||z-\tilde{z}||_{L^1(\Omega)} ~dt ,
\end{split}
\end{equation}
where 
$$
A := \frac{r(1-\delta)}{\displaystyle \int_\Omega u} - \frac{\chi^2 a^2}{8D} +  \frac{a}{||u||_{L^\infty(\Omega \times (0,\infty))}} >0,
$$
due to the choice of $\delta$ in \eqref{delta}. Taking into account that $\displaystyle \int_\Omega u$ and $\tilde{u}$ are bounded from below by Lemmas \ref{l3.1}, and \ref{l3.2}, $F_2(\rho) \geq0$, and $-F_3(\rho) \geq -c_3$ for all $\rho>0$ as previously noted, to guarantee the integrability of $k_1$, we only have to check that
$$ \int_0^\infty ||z-\tilde{z}||_{L^1(\Omega)} ~dt + \int_0^\infty \int_\Omega \left(f - \int_\Omega f\right)^2 <\infty.$$
The first one is a direct consequence of the fact that
$$
\frac{d}{dt} \int_\Omega (z -\tilde{z}) = - \int_\Omega (z - \tilde{z}) + \int_\Omega (f-\tilde{f}),
$$
and the first assumption in \eqref{hipk3} ensuring $\displaystyle \int_0^\infty ||f-\tilde{f}||_{L^1(\Omega)} ~dt < \infty$, while for the second one
$$
\int_0^\infty \Bigg|\Bigg|f - \int_\Omega f\Bigg|\Bigg|_{L^2(\Omega)}^2 \leq 2 \int_0^\infty ||f - \tilde{f}||^2_{L^2(\Omega)} \hspace{0.1 cm} dt +2 \int_0^\infty \Big|\Big| \int_\Omega f - \tilde{f}\Big|\Big|^2_{L^2(\Omega)}  \hspace{0.1 cm} dt < \infty,
$$
as a consequence of the triangle inequality for the $L^2$ norm and the second part of hypothesis \eqref{hipk3}. Thus, taking the limit as $\rho \to \infty$ in \eqref{3.32} proves the result.
\end{proof}
Next, we prove another integrability property, this time for $\displaystyle \int_\Omega \Big(|\Delta v|^2 + |\nabla v|^2\Big)$.
\begin{lemma}\label{l3.4}
    Let $f$ be such that hypotheses \eqref{1.2} and \eqref{1.2.1} hold for a certain $\tilde{f}$. Then, if $r>0$ is large enough in the sense of the previous lemmas, there exists $c>0$ such that
    $$
    \int_0^\infty \int_\Omega \Big(|\Delta v|^2 + |\nabla v|^2\Big) dx \hspace{0.1 cm} dt \leq c <\infty.
    $$
\end{lemma}
\begin{proof}
    We begin by noting that the second equation in \eqref{1.1} can be rewritten as
    \begin{equation}\label{4.1}
        \frac{\partial v}{\partial t} - \Delta v + v -f + \int_\Omega f - \int_\Omega v = au - a \int_\Omega u + \int_\Omega \frac{\partial v}{\partial t},
    \end{equation}
    for all $t>0$, due to having
    $$
   \int_\Omega \frac{\partial v}{\partial t} =  a \int_\Omega u  - \int_\Omega v + \int_\Omega f,
    $$
    as a result of directly integrating the equation over $\Omega$. Thus, multiplying equation \eqref{4.1} by $-\Delta v$ and integrating over $\Omega$, we obtain
    \begin{equation}\label{4.2}
    \begin{split}
       -\int_\Omega & \left[ \left(\frac{\partial v}{\partial t} - \int_\Omega \frac{\partial v}{\partial t} \right) \Delta v \right] + \int_\Omega |\Delta v|^2 -\int_\Omega \left[\left(v - \int_\Omega v\right) \Delta v \right]
       \\ &= - \int_\Omega \left[\left(f -\int_\Omega f \right) \Delta v \right]-a \int_\Omega \left[\left(u - \int_\Omega u\right)  \Delta v\right], \quad \text{for all } t >0.
       \end{split}
    \end{equation}
    Next, we consider the following identities that hold for all $t>0$
    \begin{equation}\label{4.3}
        \begin{split}
            - \int_\Omega \left[ \left(\frac{\partial v}{\partial t} - \int_\Omega \frac{\partial v}{\partial t} \right) \Delta v \right] &= \frac{1}{2} \frac{d}{dt} \int_\Omega |\nabla v |^2, \\\\
            - \int_\Omega \left[\left( v - \int_\Omega v \right)\Delta v  \right] &= \int_\Omega |\nabla v|^2, \\\\
            - a\int_\Omega \left[\left( u - \int_\Omega u \right)\Delta v  \right] & \leq \frac{1}{4} \int_\Omega |\Delta v|^2 + a^2 \int_\Omega \left( u - \int_\Omega u\right)^2, \\\\
            - \int_\Omega \left[\left( f - \int_\Omega f \right)\Delta v  \right] & \leq \frac{1}{4} \int_\Omega |\Delta v|^2 + \int_\Omega \left( f - \int_\Omega f\right)^2.
        \end{split}
    \end{equation}
    The first two are a direct consequence of having $\displaystyle \int_\Omega \Delta v = 0$ due to the Neumann homogeneous boundary conditions and a standard integration by parts, while the third and fourth one follow directly from the pointwise identity $\displaystyle - AB \leq \frac{1}{4} A^2 + B^2$ for any $A,B \in \mathbb{R}$.
    \\\\
    In this way, \eqref{4.2} can be expressed as
    \begin{equation}\label{4.4}
    \begin{split}
        \frac{1}{2} \frac{d}{dt} &\int_\Omega |\nabla v |^2 + \int_\Omega |\Delta v|^2 + \int_\Omega |\nabla v|^2 
        \\
        &\leq \frac{1}{2} \int_\Omega |\Delta v|^2 +  a^2 \int_\Omega \left( u - \int_\Omega u\right)^2 + \int_\Omega \left( f - \int_\Omega f\right)^2,
    \end{split}
    \end{equation}
    for all $t>0$. Considering that $ \displaystyle \int_\Omega \left( u - \int_\Omega u \hspace{0.1cm} dx\right)^2 dx$ is precisely the function $k_1$ defined in \eqref{3.13}, whose integrability was proven in Lemma \ref{l3.3}, and the fact that $ \displaystyle \int_0^\infty\int_\Omega \left( f - \int_\Omega f \hspace{0.1cm} dx\right)^2 dx\hspace{0.1cm} dt < \infty$ due to \eqref{1.2.1}, as seen in Lemma \ref{l3.3}, integrating \eqref{4.4} in time over $(0,\infty)$ yields the desired result.
\end{proof}
The next lemma provides another integrability result, this time for $\displaystyle \int_\Omega |\nabla u|^2 \hspace{0.1 cm}dx$.
\begin{lemma}\label{l3.5}
     Assume that $f$ satisfies \eqref{1.2} and \eqref{1.2.1} for a certain $\tilde{f}$. Then, if $r>0$ is large enough in the sense of the previous lemmas, there exists $c>0$ such that
    $$
    \int_0^\infty \int_\Omega |\nabla u|^2 dx \hspace{0.1 cm} dt \leq c < \infty.
    $$
\end{lemma}
\begin{proof}
Considering functionals $F_1$ and $F_2$ as in \eqref{3.24}, but this time directly with $v$ instead of decomposing it into $w$ and $z$, we obtain after an integration by parts
\begin{equation}\label{4.5}
    \begin{split}
        \frac{d}{dt}(F_1+F_2) &= \frac{d}{dt} \left(\int_\Omega\frac{u}{\tilde{u}}\right)
        - \frac{d}{dt}\left( \int_\Omega \ln u \right) + \frac{d}{dt} \left( \ln \int_\Omega u \right) 
        \\& \leq \frac{d}{dt} \left(\int_\Omega\frac{u}{\tilde{u}}\right) + \int_\Omega \left[-\frac{|\nabla u|^2}{u^2} + \chi \frac{\nabla u \nabla v}{u} \right] 
        + F_3 \\
        & \leq \frac{d}{dt} \left(\int_\Omega\frac{u}{\tilde{u}}\right) - \frac{1}{2}\int_\Omega \frac{|\nabla u|^2}{u^2} + \frac{\chi^2}{2} \int_\Omega |\nabla v|^2 + F_3, 
    \end{split}
\end{equation}
for all $t>0$, where
$$
F_3(t) := - \frac{\displaystyle \int_\Omega u v }{\displaystyle \int_\Omega u} + \int_\Omega v + \frac{\displaystyle r \int_\Omega u(1-u)}{\displaystyle \int_\Omega u} - r \int_\Omega (1-u),
$$
having used Young's inequality for the second step in \eqref{4.5}. It is direct to check that $\displaystyle \int_0^\infty F_3(t) \hspace{0.1 cm} dt < \infty$ using the same steps as in the last part of the proof Lemma \ref{l3.3}. Firstly, for all $t>0$ we have
\begin{equation*}
\begin{split}
   \int_\Omega v - \frac{\displaystyle \int_\Omega u v }{\displaystyle \int_\Omega u} &= \frac{\displaystyle \left(a \int_\Omega u + \int_\Omega f - \frac{d}{dt} \int_\Omega v \right) \cdot \left( \int_\Omega u\right) - \int_\Omega uv}{\displaystyle \int_\Omega u}   \\
   & = - \frac{d}{dt} \int_\Omega v + \frac{\displaystyle a\left( \int_\Omega u^2 -k_1 \right) + \int_\Omega u \cdot \int_\Omega f - \int_\Omega uv}{\displaystyle \int_\Omega u} \\
   & \leq - \frac{a}{||u||_{L^\infty(\Omega \times (0,\infty))}} \cdot k_1 + \left( \frac{||u||_{L^\infty(\Omega \times (0,\infty))}}{c_1} - 1 \right) \frac{d}{dt} \int_\Omega v, 
\end{split}
\end{equation*}
where $c_1>0$ is such that $c_1 \leq \displaystyle \int_\Omega u(x,t) \hspace{0.1 cm} dx$ for all $t>0$, as provided by Lemma \ref{l3.1}. 
\\\\Moreover, for the second part of $F_3$, we have
\begin{equation*}
\begin{split}
   \frac{\displaystyle r \int_\Omega u(1-u)}{\displaystyle \int_\Omega u} - r \int_\Omega (1-u) = r \int_\Omega u - \frac{\displaystyle r \int_\Omega u^2}{\displaystyle \int_\Omega u} = - \frac{r k_1}{\displaystyle \int_\Omega u} \leq 0 , \quad \text{for all } t>0.
\end{split}
\end{equation*}
Thus, it follows that $\displaystyle \int_0^\infty F_3(t) \hspace{0.1 cm} dt < \infty$. Next, a time integration of \eqref{4.5} provides
$$
\int_0^\infty \int_\Omega \frac{|\nabla u|^2}{u^2} dx \hspace{0.1 cm} dt < \infty,
$$
as $\displaystyle \int_0^\infty \int_\Omega |\nabla v|^2 dx \hspace{0.1 cm} dt$ is finite by Lemma \ref{l3.4}. The global boundedness of $u$ established in Theorem \ref{t-exst} allows us to bound our desired integral in the following way
\begin{equation}\label{4.6}
    \int_0^\infty \int_\Omega |\nabla u|^2 dx \hspace{0.1 cm} dt  \leq ||u||_{L^\infty(\Omega \times (0,\infty))}^2  \int_0^\infty \int_\Omega |\nabla v|^2 dx \hspace{0.1 cm} dt < \infty,
\end{equation}
which finishes the proof.
\end{proof}
Once the integrability of $\displaystyle \int_\Omega |\nabla u|^2 $ has been achieved, by continue by proving a pointwise bound for it.
\begin{lemma}\label{l3.6}
    Let $f$ be such that hypotheses \eqref{1.2} and \eqref{1.2.1} hold for a certain $\tilde{f}$. Then, if $r>0$ is large enough in the sense of the previous lemmas, there exists $c>0$ independent of $t$ such that
    $$
    \int_\Omega |\nabla u|^2 dx \leq c < \infty, \quad \text{for all } t>0.
    $$
\end{lemma}
\begin{proof}
We begin by multiplying the first equation in \eqref{1.1} by $-\Delta u$, taking into account that $\nabla \cdot(u \nabla v) = \nabla u \cdot \nabla v+ u \Delta v$. This yields
\begin{equation}\label{4.7}
\begin{split}
    \frac{d}{dt}\frac{1}{2}\int_\Omega |\nabla u|^2  + D\int_\Omega |\nabla u|^2 = -\chi \int_\Omega \Big(\nabla u \cdot \nabla v + u \Delta v \Big) \Delta u 
    \\  - r \int_\Omega u\left(1-u-\frac{v}{r}\right)\Delta u  , \quad \text{for all } t>0,
\end{split}
\end{equation}
where we used the fact that
$$
\int_\Omega u_t \Delta u = - \frac{d}{dt}\frac{1}{2} \int_\Omega |\nabla v|^2, \quad \text{for all } t >0. 
$$
By Theorem \ref{t-exst}, $u$, $v$ and $|\nabla v|$ are uniformly bounded in $L^\infty\big(\Omega \times (0, \infty) \big)$, so combining their bounds with Young's inequality, we obtain the following estimate for the first term on right hand side of \eqref{4.7}
\begin{equation}\label{4.8}
    -\chi \int_\Omega \Big(\nabla u \cdot \nabla v + u \Delta v \Big) \Delta u  \leq \frac{D}{2} \int_\Omega |\Delta u|^2 + c_1 \int_\Omega |\nabla u|^2+ c_2 \int_\Omega |\Delta v|^2, 
\end{equation}
for all $t>0$, with $\displaystyle c_1 := \frac{\chi^2}{D} ||\nabla v||_{L^\infty(\Omega \times (0, \infty))}>0$, and $c_2 := \displaystyle\frac{\chi^2}{D} ||u||_{L^\infty(\Omega \times (0, \infty))}>0$.
\\\\
Similarly, integrating by parts and using Young's inequality, we obtain
\begin{equation}\label{4.9}
    \begin{split}
      - r \int_\Omega u\left(1-u-\frac{v}{r}\right)\Delta u &= r \int_\Omega (1-2u) |\nabla u|^2  + \int_\Omega v |\nabla u|^2 + \int_\Omega \nabla u \cdot \nabla v
      \\  &\leq c_3 \int_\Omega |\nabla u|^2 + c_4 \int_\Omega |\nabla v|^2,   \quad \text{for all } t >0,    
    \end{split}
\end{equation}
where $c_3 := \displaystyle r \left(\frac{3}{2}+2 ||u||_{L^\infty(\Omega \times (0, \infty))} + ||v||_{L^\infty(\Omega \times (0, \infty))}\right) >0$, $c_4 = \displaystyle \frac{1}{2} >0$.
\\\\
Thus, substituting bounds \eqref{4.8} and \eqref{4.9} into \eqref{4.7}, we obtain
$$
\frac{d}{dt} \frac{1}{2} \int_\Omega |\nabla u|^2 + \frac{D}{2} \int_\Omega |\Delta u|^2 \leq (c_1 + c_3) \int_\Omega |\nabla u|^2 + c_2 \int_\Omega |\Delta v|^2 + c_4 \int_\Omega |\nabla v|^2,
$$
for all $t>0$, after which a direct integration provides the result, taking into account that the all the terms appearing on the right hand side are integrable by means of Lemmas \ref{l3.4} and \ref{l3.5}. Notice that as a by-product, we also obtain that
$$
\int_0^\infty \int_\Omega |\Delta u|^2 dx \hspace{0.1 cm} dt < \infty.
$$
\end{proof}
Next, after the integrability properties obtained in Lemmas \ref{l3.4}-\ref{l3.6}, we can prove a key result concerning the boundedness of the time derivative of the function $k_1$ introduced in \eqref{3.13}.
\begin{lemma}\label{l3.7}
    Let $k_1$ be the function defined in \eqref{3.13} and assume that $f$    
    satisfies \eqref{1.2} and \eqref{1.2.1} for a certain $\tilde{f}$. If $r>0$ is large enough in the sense of the previous lemmas, then there exists $c>0$ independent of $t$ such that
    $$
    \left| \frac{d}{dt}k_1(t)\right| \leq c < \infty, \quad \text{for all } t>0.
    $$
\end{lemma}
\begin{proof}
In view of the regularity properties of the classical solution $(u,v)$ assessed in Theorem \ref{t-exst}, we have that $k_1 \in C^1(0,\infty)$, and computing its derivative, we obtain
\begin{equation}\label{4.10}
\begin{split}
    \frac{1}{2} \frac{d}{dt} k_1(t) &= \frac{1}{2} \frac{d}{dt} \int_\Omega \left( u - \int_\Omega u\right)^2 = \int_\Omega  \left( u - \int_\Omega u \right) u_t 
    \\
    & =  D \int_\Omega \left[ \left( u - \int_\Omega u\right)  \Delta u \right] + \chi \int_\Omega \left[ \left( u - \int_\Omega u\right) \nabla \cdot(u \nabla v) \right]      
    \\
    & + r \int_\Omega \left[ \left( u - \int_\Omega u\right)   u \left(1-u-\frac{v}{r} \right) \right] , \quad \text{for all } t >0.
\end{split}
\end{equation}
Next, we compute the absolute value of the terms appearing in \eqref{4.10}. Integrating by parts and due to the Neumann homogeneous boundary conditions we have
$$
\left |D \int_\Omega \left[ \left( u - \int_\Omega u\right)  \Delta u \right] \right|= \left |-D \int_\Omega |\nabla u|^2 - D \int_\Omega u \cdot \int_\Omega \Delta u \right| =  D \int_\Omega |\nabla u|^2,
$$
as well as 
\begin{equation*}
\begin{split}
  \left|\hspace{0.1 cm} \chi \int_\Omega \left[ \left( u - \int_\Omega u\right) \nabla \cdot(u \nabla v) \right] \right|  = \left| - \chi \int_\Omega u \nabla u \cdot \nabla v - \chi \int_\Omega u \cdot \int_\Omega \nabla \cdot (u\nabla v) \right|
   \\
    = \left|- \chi \int_\Omega u \nabla u \cdot \nabla v \right| \leq \frac{\chi^2}{2} ||\nabla v||_{L^\infty(\Omega \times (0,\infty))}^2 + \frac{\chi^2}{2}||u||_{L^\infty(\Omega \times (0,\infty))}^2 \int_\Omega |\nabla u|^2, 
\end{split}
\end{equation*}
for all $t>0$, where we used Young's inequality and the $L^\infty\big(\Omega \times (0,\infty)\big)$ bounds for $u$ and $|\nabla v|$ for the last step.
\\\\
Similarly, we have
\begin{equation*}
    \begin{split}
        \left |\hspace{0.1 cm} r \int_\Omega \left[ \left( u - \int_\Omega u\right)  u \left(1-u-\frac{v}{r} \right) \right]  \right| \leq r \int_\Omega u^2 \left |1-u - \frac{v}{r} \right| + r \int_\Omega u \cdot \int_\Omega u \left|1-u - \frac{v}{r} \right|
        \\
        \leq 2r ||u||_{L^\infty(\Omega \times (0,\infty))}^2 \cdot\left(1 +||u||_{L^\infty(\Omega \times (0,\infty))} \right) + 2||u||_{L^\infty(\Omega \times (0,\infty))}^2 \cdot ||v||_{L^\infty(\Omega \times (0,\infty))},
    \end{split}
\end{equation*}
for all $t>0$. In this way, a combination of the three estimates of the right hand side of \eqref{4.10}, and fact that by Lemma \ref{l3.6}, $\displaystyle \int_\Omega |\nabla u|^2$ is bounded for all $t>0$, leads to the boundedness of $\displaystyle\left| \frac{d}{dt}k_1(t)\right|$ for all $t>0$.
\end{proof}
\section{Asymptotic behavior: convergence to the associated ODE system}\label{s4}
The main result in this work concerns the asymptotic behavior of the solutions to the original system \eqref{1.1}. In particular, the convergence to the solution associated ODE system \eqref{1.3}, as stated in the following theorem.
\begin{theorem}\label{t-per}
    Let $\Omega \subset \mathbb{R}^n$ be a bounded convex domain with smooth boundary and $f$ be a positive function satisfying \eqref{1.2} and \eqref{1.2.1}. Then, for all $\chi \in \mathbb{R}$, $D>0$ and $a\geq 0$ there exists $r_0 >0$ such that if $r\geq r_0$, system \eqref{1.1} under initial and boundary conditions \eqref{1.1.2} admits a unique solution $(u,v)$ such that
    $$
    ||u-\tilde{u}||_{L^2(\Omega)} + ||v-\tilde{v}||_{L^2(\Omega)} \to 0, \quad \text{as } t \to \infty,
    $$
    where $(\tilde{u},\tilde{v})$ is the solution to the ODE system \eqref{1.3} with initial values \eqref{1.4}.
\end{theorem}
The first stage of the proof is to show that, using the results from the previous section, $u$, converges to $\bar{u}$, given by 
$$
\bar{u}(t) := \int_\Omega u(x,t) \hspace{0.1 cm} dx, \quad  \text{for all } t>0.
$$
Subsequently, the second part consists in proving that $\bar{u}$ indeed converges to $\tilde{u}$, the solution to the ODE system \eqref{1.3}, and thus $u$ will converge to $\tilde{u}$. A similar argument will allow us to directly obtain the convergence of $v$ to $\tilde{v}$. 
\\\\
The following Lemma plays a key role in assessing the convergence.
\begin{lemma}[Lemma 5.1 in Friedman-Tello \cite{Friedman-Tello}]\label{F-T} Let $k: [0,\infty) \to \mathbb{R}$ be a $C^1$ function satisfying
\begin{enumerate}
    \item $k(t) \geq 0$ for all $t>0$,
    \item There exists $C>0$ such that $|k'(t)| \leq C$ for all $t>0$,
    \item $\displaystyle \int_0^\infty k(t) \hspace{0.1 cm} dt < \infty$,
\end{enumerate}
then $k(t) \to 0$ as $t \to \infty$.    
\end{lemma}
With it, we finish the section with the proof of Theorem \ref{t-per}.
\\\\
\textit{Proof of Theorem \ref{t-per}.} We first apply Lemma \ref{F-T} to the function
$$  k_1(t) := \int_\Omega \left(u(x,t) - \int_\Omega u(x,t)~ dx  \right)^2dx,
$$
defined in \eqref{3.13}. Clearly, $k_1(t)\geq0$ for any $t>0$, its time derivative is bounded by Lemma \ref{l3.7} and the function is integrable in time by Lemma \ref{l3.3}. Thus, 
$$k_1(t) := \Bigg|\Bigg| \displaystyle u-\int_\Omega u \Bigg|\Bigg| _{L^2(\Omega)}^2 \to 0, \quad \text{as } t \to \infty.$$
Next, to prove that $\displaystyle \int_\Omega u(x,t) \hspace{0.1 cm} dx$ converges to $\tilde{u}$, we define
$$
k_2(t) := \left(\int_\Omega u(x,t) \hspace{0.1 cm} dx - \tilde{u} \right)^2,
$$
and consider $F_2(t)$ as given by \eqref{3.24}, this is
$$
F_2(t) := \ln \left( \int_\Omega u(x,t) \hspace{0.1 cm} dx \right) - \ln \tilde{u}(t), \quad \text{for all } t>0.
$$
Moreover, we also introduce 
\begin{equation}\label{4.10.2}
    k_3(t) := \int_\Omega \big(v(x,t)-\tilde{v}(t)\big)^2 \hspace{0.1 cm} dx.
\end{equation}
Our aim is to apply Lemma \ref{F-T} again to $k_2$ and $k_3$ and thus grant the desired convergence.
\\\\
We begin by computing the time derivative of $F_2$, this time without splitting $v$ into $w$ and $z$ as done in Lemma \ref{l3.3}. We obtain
\begin{equation*}
\begin{split}
    \frac{d}{dt}F_2(t) &= \frac{\displaystyle \int_\Omega \frac{\partial u}{\partial t}}{\displaystyle \int_\Omega u} - \frac{\displaystyle \frac{d \tilde{u}}{dt~}}{\tilde{u}} =  \frac{\displaystyle-\int_\Omega \left(ru^2 -uv\right) }{\displaystyle \int_\Omega u} + r\tilde{u}+\tilde{v} 
    \\
    & = - \frac{\displaystyle r \int_\Omega u(u-\tilde{u})}{\displaystyle\int_\Omega u}-\frac{\displaystyle\int_\Omega u (v-\tilde{v})}{\displaystyle \int_\Omega u} , \quad \text{for all } t>0,
\end{split}
\end{equation*}
and thus,
\begin{equation}\begin{split}\label{4.11}
    \frac{d}{dt}F_2 + \frac{\hspace{0.3 cm}r \tilde{u}}{\displaystyle \int_\Omega u}\left(\int_\Omega u - \tilde{u} \right) &= - \frac{\displaystyle \int_\Omega u (v-\tilde{v})}{\displaystyle \int_\Omega u} + \frac{\hspace{0.1 cm} r}{\displaystyle \int_\Omega u} \left( -\int_\Omega u^2 +2\tilde{u}\int_\Omega u -(\tilde{u})^2 \right)
    \\
    & \leq -\frac{\displaystyle \int_\Omega u (v-\tilde{v})}{\displaystyle \int_\Omega u} + \frac{r k_1}{\displaystyle \int_\Omega u}, \quad \text{for all } t>0,
\end{split}    
\end{equation}
where we used Young's inequality and the fact that
\begin{equation*}
\begin{split}
    -\int_\Omega u^2 +2\tilde{u}\int_\Omega u -(\tilde{u})^2 = - \int_\Omega u^2 - \left(\int_\Omega u - \tilde{u}\right)^2 + \int_\Omega u^2 
    \\
     \leq  - \int_\Omega u^2 + \left( \int_\Omega u \right)^2 = k_1, \quad \text{for all } t>0.
    \end{split}
\end{equation*}
Multiplying \eqref{4.11} by $F_2$, using the global boundedness of $u$ and the lower bound for $\displaystyle \int_\Omega u$ proved in Lemma \ref{l3.1}, we obtain
\begin{equation}\label{4.11.2}
\begin{split}
  \frac{1}{2} \frac{d}{dt} F_2 + \frac{\hspace{0.3 cm}r \tilde{u}}{\displaystyle \int_\Omega u}\left(\int_\Omega u - \tilde{u} \right) F_2 \leq  \left( -\frac{\displaystyle \int_\Omega u (v-\tilde{v})}{\displaystyle \int_\Omega u} + c_1 k_1 \right)F_2 
  \\\leq  -\frac{\displaystyle \int_\Omega u (v-\tilde{v})}{\displaystyle \int_\Omega u} F_2 + c_2 k_1,  \quad \text{for all } t>0.
\end{split}
\end{equation}
for certain positive constants $c_1, c_2>0$. The last step follows from the fact that $F_2$ is bounded from above as a consequence of the global boundedness of $u$ and the lower bound for $\tilde{u}$ provided by Lemma \ref{l3.2}.
\\\\
Next, the Mean Value Theorem applied to the function $\ln(s)$ provides the following identity
$$
\ln(b) - \ln(a) = \frac{1}{c} (b-a),  
$$
for any $0<a<b$ and a certain $c \in (a,b)$. Thus, as $F_2 = \displaystyle \ln \left( \int_\Omega u \right) - \ln \tilde{u}$, using the identity and multiplying the result again by $F_2$ yields
\begin{equation} \label{vm-1}
    \xi(t) F_2^2 = \left(\int_\Omega u - \tilde{u}\right) F_2, \quad \text{for all } t>0,
\end{equation}
for a certain $\xi(t)$ that lies between the values of $\tilde{u}(t)$ and $\displaystyle \int_\Omega u(x,t) \hspace{0.1 cm} dx$, this is 
\begin{equation}\label{vm-2}
    \xi(t) \in \left(\min \left\{ \tilde{u}(t), \int_\Omega u(x,t) \hspace{0.1 cm} dx \right\}, ~\max \left\{ \tilde{u}(t), \int_\Omega u(x,t) \hspace{0.1 cm} dx \right\}  \right), 
\end{equation}
for all $t>0$. As both $\tilde{u}(t)$ and $\displaystyle \int_\Omega u(x,t) \hspace{0.1 cm} dx $ are globally bounded from below by a positive constant due to Lemmas \ref{l3.1} and \ref{l3.2}, there exists $c_3>0$ independent of $t$ such that
\begin{equation} \label{4.12.0}
    c_3 F_2^2 \leq \left(\int_\Omega u - \tilde{u}\right) F_2, \quad \text{for all } t>0.
\end{equation}
Moreover, the factor $\displaystyle \frac{\hspace{0.3 cm}\tilde{u}}{\displaystyle \int_\Omega u}$ is also uniformly bounded from below, so \eqref{4.11.2} can be rewritten as
\begin{equation} \label{4.12}    
 \frac{1}{2} \frac{d}{dt} F_2^2 + c_4 \cdot r \cdot F_2^2 \leq -\frac{\displaystyle \int_\Omega u (v-\tilde{v})}{\displaystyle \int_\Omega u} F_2 + c_2 k_1,  \quad \text{for all } t>0,
\end{equation}
for a certain constant $c_4>0$. Next, we bound the first term on the right hand side of the above expression. By denoting by $\varepsilon_1>0$ the lower bound for $\displaystyle \int_\Omega u $ provided by Lemma \ref{l3.1}, the Cauchy-Schwarz inequality combined with Young's inequality yields
\begin{equation*}
    \begin{split}
      -\frac{\displaystyle \int_\Omega u (v-\tilde{v})}{\displaystyle \int_\Omega u} F_2 & \leq \frac{1}{\varepsilon_1} \cdot |F_2| \int_\Omega u |v-\tilde{v}| \leq \frac{1}{\varepsilon_1} \cdot |F_2| \cdot ||u||_{L^2(\Omega)} \cdot ||v-\tilde{v}||_{L^2(\Omega)}  \\
      & \leq \frac{1}{\varepsilon_1} \left( \frac{\varepsilon_1}{2} \cdot ||v-\tilde{v}||_{L^2(\Omega)}^2 + \frac{1}{2 \varepsilon_1} \cdot ||u||_{L^2(\Omega)}^2 \cdot F_2^2 \right) \\\\
     & \leq \frac{1}{2} \int_\Omega (v-\tilde{v})^2 + c_5 F_2^2, \quad \text{for all } t>0,
    \end{split}
\end{equation*}
for $c_5 := \displaystyle \frac{1}{2 \varepsilon_1^2} \cdot  ||u||_{L^\infty(\Omega \times (0, \infty))}^2 >0$. Thus, \eqref{4.12} now becomes
\begin{equation}\label{4.13}
     \frac{1}{2} \frac{d}{dt} F_2^2 + (c_4 \cdot r - c_5) \cdot F_2^2 \leq   c_2 k_1 + \frac{1}{2} \int_\Omega (v-\tilde{v})^2,  \quad \text{for all } t>0.
\end{equation}
The next step requires bounding $\displaystyle \int_\Omega (v-\tilde{v})^2$. To do so, by subtracting the equations for $v$ and $\tilde{v}$, we obtain
$$
\frac{d}{dt} (v-\tilde{v}) + (v-\tilde{v})= \Delta v + a(u-\tilde{u}) + (f-\tilde{f})  , \quad \text{for all } t>0,
$$
and taking squares on both sides and integrating over $\Omega$ yields
\begin{equation}\label{4.14}
    \begin{split}
        \int_\Omega &\left(\frac{d}{dt} (v-\tilde{v}) \right)^2 + \int_\Omega (v-\tilde{v})^2 + \frac{d}{dt} \int_\Omega (v-\tilde{v})^2 \\
        &\leq 3 a^2 \int_\Omega (u-\tilde{u})^2 + 3 \int_\Omega |\Delta v|^2 + 3 \int_\Omega (f-\tilde{f})^2, \quad \text{for all } t >0. 
    \end{split}
\end{equation}
Thus, adding together \eqref{4.13} and \eqref{4.14} results in
\begin{equation}\label{4.15}
\begin{split}
   \frac{1}{2} \frac{d}{dt}& F_2^2 + (c_4 \cdot r - c_5) \cdot F_2^2 +   \int_\Omega \left(\frac{d}{dt} (v-\tilde{v}) \right)^2 + \frac{1}{2} k_3 + \frac{d}{dt} k_2   
   \\ &   \leq 3 a^2 \int_\Omega (u-\tilde{u})^2 + 3 \int_\Omega |\Delta v|^2 + 3 \int_\Omega (f-\tilde{f})^2 + c_2 k_1,
\end{split}
\end{equation}
for all $t>0$, where recall that we had defined $k_3 := \displaystyle \int_\Omega (v-\tilde{v})^2$ in \eqref{4.10.2}. The final part consists in estimating $\displaystyle\int_\Omega (u-\tilde{u})^2$. To this end, we have
\begin{equation}\label{4.16}
\begin{split}
    &\int_\Omega (u-\tilde{u})^2  = \int_\Omega \left[\left(u-\int_\Omega u\right) + \left(\int_\Omega u - \tilde{u} \right) \right]^2 
    \\\\
    &= \int_\Omega \left(u-\int_\Omega u\right)^2 + \int_\Omega \left(\int_\Omega u - \tilde{u} \right)^2 + 2 \int_\Omega \left(u-\int_\Omega u\right)\left(\int_\Omega u - \tilde{u} \right) 
    \\\\
    & = \int_\Omega \left(u-\int_\Omega u\right)^2  + \left(\int_\Omega u - \tilde{u} \right)^2 = k_1 + k_2, \quad \text{for all } t>0.
\end{split}
\end{equation}
Notice as well due to \eqref{vm-1}, bounding the resulting $\xi(t)$ in \eqref{vm-2} from above (is is direct to check that $\tilde{u}(t) \leq \max\{1, \tilde{u}(0)\}$ for all $t>0$), there exists $C >0$ such that
$k_2 \leq C F_2^2$. In this way, after substituting in \eqref{4.15} we obtain
\begin{equation}\label{4.15}
\begin{split}
   \frac{1}{2} \frac{d}{dt}& F_2^2 + \big(C(c_4 \cdot r - c_5) - 3a\big) \cdot k_2 +   \int_\Omega \left(\frac{d}{dt} (v-\tilde{v}) \right)^2 + \frac{1}{2} k_3 + \frac{d}{dt} k_2   
   \\ &   \leq (3 a^2 +c_2)k_1 + 3 \int_\Omega |\Delta v|^2 + 3 \int_\Omega (f-\tilde{f})^2, \quad \text{for all } t >0.
\end{split}
\end{equation}
Thus, for a given $a>0$, if $r$ is large enough so that $C(c_4 \cdot r - c_5) - 3a >0$, a direct time integration results in
$$
\int_0^\infty k_2(t) \hspace{0.1 cm} dt + \int_0^\infty k_3(t) \hspace{0.1 cm} dt < \infty,
$$
as a consequence of the integrability of all the terms on the right hand side of \eqref{4.15} due to Lemma \ref{l3.3}, Lemma \ref{l3.4} and hypothesis \eqref{1.2.1}.
\\\\
Lastly, it is direct to check that $\displaystyle \left|\frac{d}{dt} k_2 \right|$ and $\displaystyle \left|\frac{d}{dt} k_3 \right|$ are indeed globally bounded, so by Lemma \ref{F-T}, $k_2(t) \to 0$ and $k_3 \to 0$ as $t \to \infty$. 
\\\\
As we had already proved that $k_1(t) \to 0$ at the beginning of the proof, $k_2(t) \to 0$ directly implies that $||u-\tilde{u}||_{L^2(\Omega)} \to 0$ as $t \to \infty$. The same result for $v$ is a direct consequence of $k_3(t) \to 0$. \qedsymbol{}

\section{Periodicity}\label{s5}
Lastly, once the convergence to $(\tilde{u},\tilde{v})$ has been established, we include some final remarks regarding its qualitative properties in the case in which $\tilde{f}$ is a periodic function. The natural question that arises is whether the solution $(\tilde{u},\tilde{v})$ does inherit the periodicity of $\tilde{f}$ or not.  
\\\\
As mentioned in the Introduction, the ODE system \eqref{1.3} was considered in \cite{HN25} for $a=0$, and a periodic function $\tilde{f}$ of period $T>0$, having
\begin{equation}\label{5.1}
    \begin{cases}
  \displaystyle  \frac{d \tilde{u}}{dt~} = r \tilde{u} (1 - \tilde{u}) - \tilde{u}\tilde{v}, \quad & t>0,\\\\
  \displaystyle  \frac{d \tilde{v}}{dt~} = - \tilde{v} + \tilde{f},\quad & t>0, 
\end{cases}
\end{equation}
with initial values
$$
   \tilde{u}(0) = \int_\Omega u_0(x) ~dx, \quad \tilde{v}(0) = \int_\Omega v_0(x) ~dx.
$$
As the second equation is uncoupled, the system can be explicitly solved, determining the threshold value 
\begin{equation} \label{5.2}
r_{\min} := \displaystyle \frac{1}{T} \int_0^T \tilde{f}(s) \hspace{0.1 cm} ds,
\end{equation}
the average value of $\tilde{f}$ over its period. It is proven in the paper that if $r>r_{\min}$, then system \eqref{5.1} admits a unique positive periodic solution, that can be constructed by adequately selecting $\tilde{u}(0)$ and $\tilde{v}(0)$. In particular, the unique periodic solution corresponds to the pair of initial values $(\tilde{u}(0), \tilde{v}(0)) = (u^{per}_0,v^{per}_0)$, where
\begin{equation}\label{5.3}
\begin{split}
    u^{per}_0(r) &= \frac{e^{\displaystyle \int_0^T (r-\tilde{f}(s)) ~ds} - 1}{ \displaystyle \int_0^T \left( r ~e^{\displaystyle \int_0^s (r- \tilde{f}(\tau))~ d \tau} \right) ~ds},  \\\\
v^{per}_0(r) &=\frac{1}{e^T - 1} \int_0^T \tilde{f}(s) e^s ~ds.
\end{split}
\end{equation}
\\\\
For $a>0$, the question remains open due to the nonlinear coupling of the system. The result obtained in Lemma \ref{l3.2} provides a partial answer to this: if $r>r_{\min, a}$, where
\begin{equation} \label{5.4}
r_{\min, a} := \max \big \{ \tilde{v}(0), ~ a \cdot \max\{\tilde{u}(0), 1\} + ||\tilde{f}||_{L^\infty (0, \infty)} \big\},
\end{equation}
then there exists a constant $c>0$ such that $\tilde{u}(t)>c$ for all $t>0$. While this rules out a decay of $\tilde{u}$ to zero, it does not directly imply that $(\tilde{u},\tilde{v})$ is $T$–periodic. Moreover, $r_{\min,a}$ is not a sharp threshold in the sense of the decoupled case, since substituting $a=0$ in $r_{\min,a}$ produces a larger value than the exact threshold $r_{\min}$ obtained for $a=0$.

\section*{Acknowledgments}
This work was supported by Project PID2022-141114NB-I00 from the Spanish Ministry of Science and Innovation (M.N) and by Grant FPU23/03170 from the Spanish Ministry of Science, Innovation and Universities (F.H.-H.) 

\section*{Declaration of Competing Interest}
The authors declare that they have no known competing financial interests or personal relationships that could have appeared to influence the work reported in this paper.

\end{document}